\documentclass[11pt,oneside,a4paper]{article}
\usepackage[T1]{fontenc}
\usepackage[utf8]{inputenc}
\usepackage{graphicx}
\usepackage{epsfig}
\usepackage{psfrag}
\usepackage[table,dvipsnames]{xcolor}
\usepackage{pst-node}
\usepackage{pstricks}
\usepackage[english]{babel}
\RequirePackage{amssymb}
\RequirePackage{amsmath}
\RequirePackage{amsfonts}
\RequirePackage{amscd}
\RequirePackage{mathrsfs}
\RequirePackage{multirow}
\newcommand{\ol}{\overline}

\newcommand{\ub}{\underbrace}
\newcommand{\ds}{\displaystyle}
\newcommand{\de}{\delta}

\newcommand{\la}{\lambda}
\renewcommand{\rho}{\varrho}
\renewcommand{\phi}{\varphi}

\newcommand{\setN}{\mathbb{N}} 
\newcommand{\setR}{\mathbb{R}} 
\newcommand{\setC}{\mathbb{C}} 
\newcommand{\setS}{\mathbb{S}} %
\newcommand{\Union}{\bigcup}

\newcommand{\Intersection}{\bigcap}

\newcommand{\dsumm}{\bigoplus\limits} 
\newcommand{\summ}{\sum\limits}

\newcommand{\limm}{\lim\limits}

\renewcommand{\dim}{\operatorname{dim}}

\newcommand{\pol}{\operatorname{P}} 

\newcommand{\minn}{\min\limits}
\newcommand{\norm}[1]{\Vert #1 \Vert}
\newcommand{\rank}{\operatorname{rank}} 
\newcommand{\e}{\operatorname{e}} 
\newcommand{\innp}[2]{\left\langle #1,#2\right\rangle} 
\newcommand{\Span}[1]{\operatorname{span}\left\{#1\right\}} 
\newcommand{\intt}[2][]{
    \int\limits #2
    \def\testempty{#1}
    \ifx\testempty\empty
        \operatorname{dx}
    \else
        \operatorname{d#1}
    \fi
}
\newcommand{\Fmla}[1]{g(#1)}
\newcommand{\WA}{\mathcal{A}}
\newcommand{\Bf}{{\bf f}}
\newcommand{\Btf}{{\bf\tilde f}}
\newcommand{\Ba}{{\bf a}}
\newcommand{\Bta}{{\bf\tilde a}}
\newcommand{\BPhi}{{\bf \Phi}}

\newcommand{\sph}{{\mathcal H}} 
\newcommand{\closure}[1]{\operatorname{closure}\left(#1\right)}
\newtheorem{DfnEnv}{Definition}[section]
\newcommand{\dfn}[1]{\begin{DfnEnv}#1\end{DfnEnv}}
\newtheorem{ThmEnv}{Theorem}[section]
\newcommand{\thm}[1]{\begin{ThmEnv}#1\end{ThmEnv}}
\newtheorem{LmmEnv}{Lemma}[section]
\newcommand{\lmm}[1]{\begin{LmmEnv}#1\end{LmmEnv}}
\newcommand{\set}[2][]{
    \def\testempty{#1}
    \ifx\testempty\empty
        \left\{#2\right\}
    \else
        \left\{#2\quad{\bf :}\quad #1\right\}
    \fi
}
\newcommand{\Set}[2][]{
    \def\testempty{#1}
    \ifx\testempty\empty
        \Big\{#2\Big\}
    \else
        \Big\{#2\quad{\bf :}\quad #1\Big\}
    \fi
}
\newlength{\textlen}
\newcommand{\bignu}{\text{{\LARGE$\nu$}}}
\newcommand{\wlsp}{ 
   \settowidth{\textlen}{\ensuremath{\bignu}}
  \ensuremath{\bignu\hspace*{-0.3\textlen}\bignu}
}
\title{Frame and wavelet systems on the sphere}
\author{{\sc Margit Pap}\\ \ \\ University of P\'ecs, Ifj\'us\'ag \'utja 6, 7634 P\'ecs, Hungary\\
papm@ttk.pte.hu}
\date{}

\begin{document}
\pagestyle{empty}
\maketitle
\bigskip
\begin{abstract}
In this paper we formulate a weighted version of  minimum problem (1.4) on the sphere and
 we show that, for $K\le L$, if $\set{\phi_k}^K_{k=1}$ consists of
the spherical functions with degree less than $N$ we can localize the points
$(\xi_1,\ldots,\xi_L)$ on the sphere so that the solution of this problem is the simplest possible. This localization is connected
to the discrete orthogonality of the spherical functions which was proved in \cite{pap}.
Using these points we construct a frame system and a wavelet system on the sphere and  we study the properties of these systems. For $K>L$ a similar construction was made in paper \cite{prestin}, but in that case the solution of the minimum problem (1.4) is not as efficient as it is in our case. The analogue of Fej\'er and de la Val\'ee-Poussin summation methods introduced in \cite{pap} can be expressed by the frame system introduced in this paper.
\end{abstract}
\bigskip
\cleardoublepage
\numberwithin{equation}{section}
\pagestyle{plain}
\section{Approximation on the sphere}\label{sec1}
Let $\setS^2$ denote the three dimensional unit sphere, let $\bignu$ be a finite dimensional subspace of $L^2(\setS^2)$ with $\dim
\bignu=K$. With respect to a basis $\set{\phi_k}^K_{k=1}$ of $\bignu$ any
$f\in\bignu$ has a unique representation
\begin{equation}\label{label1}
f=\summ_{k=1}^Ka_k\phi_k,\qquad a_k\in\setC.
\end{equation}
For $\xi_l\in\setS^2$, $l=1,\ldots,L$ let suppose that we know the
values of $f(\xi_l)$, $l=1,\ldots,L$.
\\
Denote by
\begin{equation}\label{label2}
\Bf:=(f(\xi_1),\ldots,f(\xi_L))^T\in\setC^L,
\end{equation}
and
\begin{equation}\label{label4}
\BPhi:=\left(\begin{matrix}%
\phi_1(\xi_1)&\ldots&\phi_K(\xi_1)\\
\vdots&&\vdots\\
\phi_1(\xi_L)&\ldots&\phi_K(\xi_L)\\
\end{matrix}\right)\in\setC^{L\times K}.
\end{equation}
\\
Let us consider the approximation problem:
\begin{equation}\label{label5}
\text{Find } \Bta\in\setC^K%
\text{ with } \norm{\Bf-\BPhi\Bta}_2\le\norm{\Bf-\BPhi \Ba}_2%
\text{ for all } \Ba\in\setC^K
\end{equation}
or, equivalently,
$\minn_{\Ba\in\setC^k}\norm{\Bf-\BPhi\Ba}_2$,
where $\norm{\circ}_2$ denotes the normal euclidian norm.
For $K\le L$ the solution of the approximation problem \eqref{label5}
 can be found using the least squares method, and $\Bta$ is
solution of the normal equations
\begin{equation}\label{label6}
\BPhi^H\BPhi\Ba=\BPhi^H\Bf.
\end{equation}
Assuming that the matrix $\BPhi$ has full rank, i.e. $\rank(\BPhi)=K$
we obtain that
\begin{equation}\label{label7}
\Ba=(\BPhi^H\BPhi)^{-1}\BPhi^H\Bf.
\end{equation}
If  number $K$ is large, then the difficulty in computations is the
determination of $(\BPhi^H\BPhi)^{-1}$, which needs a great number of
operations.
\paragraph{Problem.}
%
The  question is: for a given system  $\set{\phi_k}_{k=1}^K$ how to
choose $(\xi_1,\ldots,\xi_L)$ so that the computation of
$(\BPhi^H\BPhi)^{-1}$  to be not difficult. The most simplest case is if
for a good choice of $(\xi_1,\ldots,\xi_L)$, the matrix
$(\BPhi^H\BPhi)^{-1}$ is equal by the identity matrix.

In what follows we will show that, for $K\le L$, if $\set{\phi_k}^K_{k=1}$ consists of
the spherical functions with degree less than $N$ we can localize the points
$(\xi_1,\ldots,\xi_L)$ in this way. This localization is connected
to the discrete orthogonality of the spherical functions.
Using these points we  will construct a frame system and a wavelet system on the sphere and  we will study the properties of these systems. For $K>L$ a similar construction was made in paper \cite{prestin}, but in that case the solution of the minimum problem (1.4) is not as efficient as it will be in our case. The analogue of Fej\'er and de la Val\'ee-Poussin summation methods introduced in \cite{pap} can be expressed by the frame system introduced in this paper.
\section{Spherical harmonics of degree $n$}\label{sec2}
Spherical harmonics play an important role in Fourier analysis. They
are restrictions to the sphere $\setS^2$ of homogeneous polynomials
that are solutions of the $3$ dimensional Laplace equation.
\\
Let denote by $\xi=(\sin \theta \cos \phi, \sin \theta \sin \phi, \cos
\theta)\in \setS^2$ and $f:\setS^2\to \setC,\, f(\xi)=f(\sin \theta
\cos \phi, \sin \theta \sin \phi, \cos \theta)$.
\\
The set of $n$ degrees spherical harmonics is denoted by $\sph_n$, with
\begin{equation}\label{label8}
\dim(\sph_n)=2n+1.
\end{equation}
For arbitrary $n\in\setN$ the functions
\begin{equation}\label{label9}
Y_{nk}(\theta,\phi):=\sqrt{2n+1}\cdot%
    \pol_n^{|k|}(\cos\theta)\e^{ik\phi},\quad k=-n,\ldots,n
\end{equation}
are called the spherical polynomials (functions) of order $n$, where $\pol_n^{|k|}$ are the associate Legendre polynomials. The spherical polynomials of order $n$
constitute an orthonormal basis for $\sph_n$, where the ortogonality
is with respect to the scalar product induced by the following
(continuous) measure on the unit sphere
\begin{equation}\label{label10}
\intt[\xi]{_{\setS^2}f(\xi)}:=\frac{1}{4\pi}%
\intt[\phi]{_0^{2\pi}\intt[\theta]{_0^\pi f(\theta,\phi)\sin\theta}},
\end{equation}
namely
$$\innp{Y_{nk}}{Y_{ml}}=\frac{1}{4\pi}%
\intt[\phi]{_0^{2\pi}\intt[\theta]{_0^\pi Y_{nk}(\theta,\phi)
    \ol{Y_{ml}(\theta,\phi)}\sin\theta}}=\de_{mn}\de_{kl}$$
(where $\de_{mn}$ is the Kronecker symbol).
\\
The addition theorem for spherical harmonics is
\begin{equation}\label{man:2:4}
\summ_{k=-n}^nY_{nk}(\xi)\ol{Y_{nk}(\eta)}=
    (2n+1)\pol_n(\xi\cdot\eta),\quad\xi,\eta\in\setS^2,
\end{equation}
 where $\pol_n$ is the Legendre polynomial of $n$--th degree.
\\
The function $K_n:\setS^2\times\setS^2\to\setR$,
\begin{equation}\label{label12}
K_n(\xi,\eta)=(2n+1)\pol_n(\xi\cdot\eta)
\end{equation}
is called reproducing kernel because it has the following
reproducing property: for arbitrary $f\in\sph_n$,
\begin{equation}\label{label13}
\innp{f}{K_n(\cdot,\eta)}=f(\eta),\text{ with }\eta\in\setS^2.
\end{equation}

\section{Discrete orthogonality of spherical functions}\label{sec3}
In paper \cite{pap} we proved that it can be constructed a set of points in
$[0,\pi]\times[0,2\pi]$ and a discrete measure, so that the
orthonormality property of spherical functions regarding to the
scalar product, induced by the discrete measure is preserved. In
what follows we summarize these results.
\\
Let denote by $\la_k^N\in(-1,1)$, $k\in\set{1,\ldots,N}$ the roots of
Legendre polynomials $\pol_N$ of order N, and for $j=1,\ldots,N$ let
\begin{equation}\label{label14}
l_j^N(x)=\frac{(x-\la_1^N)\ldots(x-\la_{j-1}^N)%
        (x-\la_{j+1}^N)\ldots(x-\la_N^N)}{%
    (\la_j^N-\la_1^N)\ldots(\la_j^N-\la_{j-1}^N)%
        (\la_j^N-\la_{j+1}^N)\ldots(\la_j^N-\la_N^N)},
\end{equation}
be the corresponding fundamental polynomials of Legendre
interpolation. Denote by
\begin{equation}\label{label15}
\WA_k^N:=\intt{_{-1}^1l_k^N(x)}\quad(1\le k\le N)
\end{equation}
the corresponding Cristoffel--numbers. Let consider the set of nodal
points
\begin{equation}\label{label16}
X:=\set[k=\ol{1,N},\ j=\ol{0,2N}]%
{z_{kj}=(\theta_k,\phi_j)=
    \left(\arccos\la_k^N,\frac{2\pi j}{2N+1}\right)}
\end{equation}
and the weights
\begin{equation}\label{label17}
\mu_N(z_{kj}):=\frac{\WA_k^N}{2(2N+1)}.
\end{equation}
These weights are positive numbers (see \cite{szego}).
On the set $X$ we consider the following discrete integral
\begin{equation}\label{label17}
\intt[\mu_N]{_Xf}:=\summ_{k=1}^N\summ_{j=0}^{2N}f(z_{kj})\mu_N(z_{kj})=
    \summ_{k=1}^N\summ_{j=0}^{2N}f(\theta_k,\phi_j)\frac{\WA_k^N}{2(2N+1)}.
\end{equation}
\thm{\label{thm1 }  {\rm (\cite{pap})} \,%
    Let $N\in\setN$, $N\ge1$, then the finite set of
    normalized spherical functions
    $$\Set[k\in\set{-n,\ldots,n},\ n\in\set{0,\ldots,N-1}]%
            {Y_{nk}:\setS^2\to\setC}$$
    form an orthonormal system on the set of modal point $X$
    regarding to the discrete integral \eqref{label17}, i.e.,
    \begin{equation}\label{man:3:6}
    \intt[\mu_N]{_XY_{nk}\ol{Y_{n'k'}}}=\de_{nn'}\de_{kk'},\quad
    n,n'<N,\ k\in\set{-n,\ldots,n},\ k'\in\set{-n',\ldots,n'}.
    \end{equation}
}
\thm{\label{thm2 } {\rm (\cite{pap})}\, %
    For all $f\in C(\setS^2)$,
    $$\limm_{N\to\infty}\intt[\mu_N]{_Xf}=
        \intt[\mu]{_{\setS^2}f}.$$
}
\section{Approximation on the sphere corresponding to points
determined by discretisation process}\label{sec4}
Let localize $\xi$ on the sphere in the following
way:
\begin{multline}\label{label19}
X'=\Big\{\xi_{kj}=(\sin\theta_k\cos\phi_j,%
        \ \sin\theta_k\sin\phi_j,%
        \ \cos\theta_k)\in\setS^2\, :\\%
    (\theta_k,\phi_j)\in X,\ k=\ol{1,N},\ j=\ol{0,2N}\Big\}.
\end{multline}
This set consists of $L=N(2N+1)$ points.
\\
Suppose that we can measure the values of $f$ on the set $X'$, namely we
know $\Bf=(f(\xi_{10}),\ldots,f(\xi_{N2N}))^T$.
\\
Let introduce the following notations:
$$I_N:=\left(\begin{matrix}%
\sqrt{\mu_N(\xi_{10})} & 0 & 0 & \dots & 0 \\
0 & \sqrt{\mu_N(\xi_{11})} & 0 & \dots & 0 \\
\vdots & & &\vdots \\
0 & 0 & 0 & \dots & \sqrt{\mu_N(\xi_{N2N})} \\
\end{matrix}\right)\in M^{L\times L},$$
$$\Bf_1=I_N\Bf=\left(\sqrt{\mu_N(\xi_{10})}f(\xi_{10}),\ \ldots,\ %
        \sqrt{\mu_N(\xi_{N2N})}f(\xi_{N2N})\right)^T,$$
$$g_{nk}(\xi)=Y_{nk}(\xi)\sqrt{\mu_N(\xi)},\quad
    n\in\{0,1,\ldots,N-1\},\quad k\in\{-n,\ldots,n\},$$
\renewcommand{\Fmla}[2]{g_{\text{\tiny{$#1$}}}(\xi_{\text{\tiny{$#2$}}})}
$$\BPhi_1=I_N\BPhi=\left(\begin{matrix}%
\Fmla{00}{10} & \Fmla{1-1}{10} & \dots & \Fmla{N-1N-1}{10} \\
\Fmla{00}{11} & \Fmla{1-1}{11} & \dots & \Fmla{N-1N-1}{11} \\
\vdots & \vdots & & \vdots \\
\Fmla{00}{N2N} & \Fmla{1-1}{N2N} & \dots & \Fmla{N-1N-1}{N2N} \\
\end{matrix}\right)\in M^{L\times N^2}.$$
Let formulate the following weighted
minimum problem:

find $\tilde a$ so that $\minn_{a\in\setC^{N^2}}\norm{\Bf_1-\BPhi_1a}=\norm{\Bf_1-\BPhi_1\tilde a}$.
 The solution  $\tilde a$, according to (1.6) is
$$\tilde a=\left((I_N\BPhi)^H(I_N\BPhi)\right)^{-1}%
        \left(I_N\BPhi\right)^H\Bf_1.$$
We show that in this case $\left((I_N\BPhi)^H(I_N\BPhi)\right)$ is
the $N^2$ dimensional identity matrix, so the computation of $\tilde
a$ is the most simplest possible. Indeed, denote
$(I_N\BPhi)^H(I_N\BPhi)=$
\renewcommand{\Fmla}[2]{A(\text{\scriptsize{$#1;\ #2$}})}
\begin{equation*}\label{label20}
\left(\begin{matrix}%
\Fmla{0,0}{0,0} & \dots & \Fmla{0,0}{N-1,N-1} \\
\vdots & & \vdots \\
\Fmla{N-1,N-1}{0,0} & \dots & \Fmla{N-1,N-1}{N-1,N-1} \\
\end{matrix}\right)\in M^{N^2\times N^2},
\end{equation*}
where
\begin{multline*}
 A(n,k;\ n',k')=
    \ol{Y_{nk}(\xi_{1,0})}Y_{n'k'}(\xi_{1,0})\mu_N(\xi_{1,0})
    +\ldots+
    \ol{Y_{nk}(\xi_{N,2N})}Y_{n'k'}(\xi_{N,2N})\mu_N(\xi_{N,2N}).
\end{multline*}
From (3.6) and the discrete orthogonality of spherical functions it follows that
\begin{multline*}
 A(n,k;\ n',k')=
    \intt[\mu_N]{_X\ol{Y_{nk}}Y_{n'k'}}=\delta_{nn'}\delta_{kk'},\\
n,n'\in\{1,\ldots,N-1\}\quad k\in\{-n,\ldots,n\},\quad%
    k'\in\{-n',\ldots,n'\},
\end{multline*}
from this it follows that

\begin{equation*}\label{label22}
(I_N\BPhi)^H(I_N\BPhi)=\left(\begin{matrix}%
1 & 0 & \dots & 0 \\
0 & 1 & \dots & 0 \\
\vdots & & \vdots \\
0 & 0 & \dots & 1 \\
\end{matrix}\right)\in M^{N^2\times N^2}.
\end{equation*}
From (1.7) and (3.5) it follows that
\begin{equation*}\label{label24}
\ds\Bta = (I_N\BPhi)^H\Bf_1=
\left(\begin{matrix}%
\ds\intt[\mu_N]{_X\ol{Y_{0 0}}f} \\
\ds\intt[\mu_N]{_X\ol{Y_{1 -1}}f} \\
\ds\vdots \\
\ds\intt[\mu_N]{_X\ol{Y_{N-1 N-1}}f} \\
\end{matrix}\right).
\end{equation*}
which means that that the components of $\Bta$ are exactly the
discrete Laplace-Fourier coefficients of $f$. In this case the best
approximant can be expressed using the reproducing kernels
$K_n(\xi,\eta)$, namely
\renewcommand{\Fmla}[1]{%
    \intt[\mu_N(\eta)]{_X\summ_{n=0}^{N-1}%
    \left(K_n(\xi_{#1},\eta)\right)%
    \sqrt{\mu_N(\xi_{#1})}f(\eta)}}
\begin{equation*}\label{label26}
\Btf=\BPhi_1\Bta=\left(\begin{matrix}%
\Fmla{1,0}\\\Fmla{1,1}\\\vdots\\\Fmla{N,2N}\\
\end{matrix}\right).
\end{equation*}
This gives the idea to study the properties of the set of functions
$$\set[\xi_{l,m}^N\in X']%
        {\sqrt{\mu_N(\xi_{l,m}^N)}\phi_j(\cdot,\xi_{l,m}^N)},$$
where
\begin{equation}\label{label41}
\phi_j(\cdot,\xi_{l,m}^N)=\summ_{n=0}^{m_j-1}\summ_{k=-n}^n\ol{Y_{n,k}(\xi_{l,m}^N)}Y_{n.k}(\cdot)=
\summ_{n=0}^{m_j-1}K_n(\cdot,\xi_{l,m}),
\end{equation}
$\{m_j\}_{j=1}^\infty$ is a strictly monotone increasing sequence
of positive integers, $j_0$ so that $m_j\le N$, if $j\le j_0$.
We will show that with this systems we can generate a multiresolution decomposition in
the space of spherical polynomials of degree $N$,
they constitute a frame system at every level of the multiresolution, and they generate a wavelet decomposition on the sphere.

\section{Multiresolution decomposition and wavelet spaces on the sphere
using the nodal points determined by $X$}
\dfn{\label{dfn1}
    A sequence $\{\bignu_j\}_{j\in\setN}$ of finite dimensional
    subspaces of $L^2(\setS^2)$ will be called a multiresolution of
    $L^2(\setS^2)$ if the following conditions are satisfied:
    \begin{equation}\label{M1}\tag{M1}
    \ds\bignu_j\subset\bignu_{j+1}\text{ for all }j\in\setN,
    \end{equation}
    \begin{equation}\label{M2}\tag{M2}
    \closure{\Union_{j\in\setN}\bignu_j,\norm{\circ}}=L^2(\setS^2).
    \end{equation}
}
Usually, a definition of multiresolution includes a condition on the
intersection of the spaces $\bignu_j$. From \eqref{M1} it follows
immediately that $$\Intersection_{j\in\setN}\bignu_j=\bignu_0.$$
We define the wavelet spaces $\wlsp_j$ as the orthogonal complements of
$\bignu_j$ in $\bignu_{j+1}$, i.e.
\begin{equation}\label{label27}
\wlsp_j:=\bignu_{j+1}\ominus\bignu_j,
\end{equation}
which means that
\begin{equation}\label{label28}
\bignu_{j+1}=\bignu_j\oplus\wlsp_j.
\end{equation}
\lmm{\label{lmm1}
    Let $\{\bignu_j\}_{j=1}^\infty$ be a multiresolution analysis of
    $L^2(\setS^2)$ and for $j\in\setN$ let $\wlsp_j$  be the
    corresponding wavelet spaces defined by \eqref{dfn1}.
    With $\wlsp_0:=\bignu_1$ we have
    \begin{equation}\label{label29}
    \ds L^2(\setS^2)=\dsumm_{j=0}^\infty\wlsp_j.
    \end{equation}
}
Defining the operators $R_j$ and $Q_j$ to be orthogonal projections
$R_j:L^2(\setS^2)\to\bignu_j$ and $Q_j:L^2\to\wlsp_j$, $j\in\setN$ we have

%

%
\begin{equation}\label{label31}
R_{j+1}=R_1f+\summ_{k=1}^jQ_kf
\end{equation}
for $f\in L^2(\setS^2)$ .
We will study an approximation process on the sphere based on harmonic
polynomials.
\thm{\label{thm3}
    Every harmonic polynomial of degree $N$ over the sphere $\setS^2$ can
    be written as a sum of spherical harmonics of degree $N$, i.e.
    \begin{equation}\label{label32}
    \Pi_N(\setS^2)=\dsumm_{n=0}^N\sph_n.
    \end{equation}
}
Particularly, any function $f\in L^2(\setS^2)$ can be approximated by
spherical harmonics in $L^2(\setS^2)$--sense up to arbitrary
precision. From $L^2(\setS^2)$--orthonormality of the space $\sph_n$
for $n\in\setN$ we obtain the dimension
\begin{equation}\label{label33}
\dim\Pi_{N-1}(\setS^2)=\summ_{n=0}^{N-1}\dim\sph_n=\summ_{n=0}^{N-1}(2n+1)=N^2.
\end{equation}
In the finite dimensional Hilbert space $\sph_n$ with the inner product
$L^2(\setS^2)$ any function $Y_n\in\sph_n$ can be represented with respect to an
orthonormal basis $\{Y_{nk}\}_{k=-n}^n$ as a Laplace--Fourier sum
\begin{equation}\label{label34}
Y_n=\summ_{k=-n}^n\innp{Y_n}{Y_{nk}}Y_{nk}.
\end{equation}
Finally, any $Y\in\Pi_{N-1}$ can be written as
\begin{equation}\label{label35}
Y(\eta)=\summ_{n=0}^{N-1}\summ_{k=-n}^n\innp{Y}{Y_{nk}}Y_{nk}(\eta).
\end{equation}
%
%
%
Because of the discrete orthonormality of spherical harmonics
 the Laplace--Fourier coefficients
$\alpha_{nk}:=\innp{Y}{Y_{nk}}$ can be computed using the discrete
integral, namely
\begin{multline}\label{man:5:10}
\alpha_{nk}=\innp{Y}{Y_{nk}}_X=\intt[\mu_N]{_X\ol{Y_{nk}}Y}=\\
\summ_{k'=1}^N\summ_{j=0}^{2N}Y(\theta_{k'},\phi_j)\cdot\ol{Y_{nk}(\theta_{k'},\phi_j)}\mu_N(z_{k'j}).
\end{multline}
This means that if we know the values of  $Y\in\Pi_{N-1}$ on the set  $X$ then we can compute the  exact values of continuous Laplace-Fourier coefficients of  $Y$, consequently we know the values of $Y$ on the whole sphere. This makes possible to develop a construction of weighted scaling functions in case $K\le L$, which is analogous to that of  J. Prestin and M. Conrad \cite{prestin}  for $K> L$.
\\
Substituting $\alpha_{n,k}$ in \eqref{label35} we obtain that
\begin{equation}\label{label38}
Y=\summ_{n=0}^{N-1}\summ_{k=-n}^n\innp{Y}{Y_{nk}}_XY_{nk},\quad
Y\in\Pi_{N-1}(\setS^2).
\end{equation}
From \eqref{label35}, \eqref{label38} and using notation
\eqref{man:2:4}, \eqref{label13} we obtain the following reproducing
formulas in the space $\Pi_{N-1}(\setS^2)$. Any
$Y\in\Pi_{N-1}(\setS^2)$ can be written in following two  ways:
\begin{equation}\label{man:5:12}
Y(\eta)=\innp{Y(\cdot)}{\summ_{n=0}^{N-1}K_n(\cdot,\eta)}=%
    \innp{Y(\cdot)}{\summ_{n=0}^{N-1}K_n(\cdot,\eta)}_X.
\end{equation}
Let $\{m_j\}_{j=1}^\infty$ be a strictly monotone increasing sequence
of positive integers. Then the spaces
$\bignu_j=\Pi_{m_j-1}(\setS^2)=\dsumm_{n=0}^{m_j-1}\sph_n$,
$(j\in\setN)$ satisfy the following conditions
\begin{enumerate}
\item $\bignu_j\subset\bignu_{j+1}$, $j\in\setN$
\item $\closure{\Union_{j=1}^\infty\bignu_j,\norm{\circ}}=L^2(\setS^2).$
\end{enumerate}
Thus $\{\bignu_j\}_{j=1}^\infty$ is a multiscale decomposition of
$L^2(\setS^2)$. We will choose  from $\bignu_1$ the scaling functions
$$\Set{Y_{nk}}_{n=0;k=-n}^{m_1-1;n},$$
and the set of nodal points $X'$ on $\setS^2$, corresponding to the
set $X$ from section 2.
\\
Let consider $j_0$ so that $m_j\le N$, if $j\le j_0$.
For the scale $\bignu_j$ we introduce the scaling functions
\begin{equation}\label{label40}
\tilde\phi_j=\summ_{n=0}^{m_j-1}\summ_{k=-n}^nY_{nk}
\end{equation}
and the ''weighted'' scaling functions from $\bignu_j$
\begin{equation}\label{label41}
\phi_j(\cdot,\xi_{l,m}^N)=\summ_{n=0}^{m_j-1}\summ_{k=-n}^n\ol{Y_{nk}(\xi_{l,m}^N)}Y_{nk}(\cdot),
\end{equation}
where $\xi_{l,m}^N\in X'$  are defined by \eqref{label19}.
\\
The advantage of this choice is the following: beside the continuous orthonormality of
the spherical functions of degree less then $N$  they have the discrete
orthonormality property on set $X$.
\\
 In this way on
the $j$--th scale the number of the  weighted scaling functions is equal to the number of the points  on which we measure the functions, namely $L=N(2N+1)$. From \eqref{man:2:4}
and \eqref{label13} it follows that
\begin{equation}\label{label42}
\phi_j(\cdot,\xi_{lm}^N)=%
\summ_{n=0}^{m_j-1}(2n+1)\pol_n(\cdot,\xi_{lm}^N)=%
\summ_{n=0}^{m_j-1}K_n(\cdot,\xi_{lm}^N)
\end{equation}
are real valued.
\\
We summarize some properties of the functions
$\phi_j(\cdot,\xi_{lm}^N)$ in the next result.
{\thm{\ \\
  \begin{enumerate}
    \item \label{item:1}
        The functions $\phi_j(\cdot,\xi_{lm}^N)$ have the
        reproducing property
        \begin{equation*}\label{label43}
        \innp{f}{\phi_j(\cdot,\xi_{lm}^N)}=f(\xi_{lm}^N)
        \text{ for all } f\in\bignu_j,\ j\le j_0
        \end{equation*}
        \begin{equation*}\label{label44}
        \innp{f}{\phi_j(\cdot,\xi_{lm}^N)}_X=f(\xi_{lm}^N).
        \end{equation*}
    \item It holds $\norm{\phi_j(\cdot,\xi_{lm}^N)}=m_j$,
        $\phi_j(\xi_{lm}^N,\xi_{lm}^N)=m_j^2$, %
        $\norm{\phi_j(\cdot,\xi_{lm}^N)}_X=m_j$.
    \item The function $\phi_j(\cdot,\xi_{lm}^N)$ is localized
        around $\xi_{lm}^N$, i.e.
        \begin{equation*}\label{label45}
        \frac{\norm{\phi_j(\cdot,\xi_{lm}^N)}}
                {\phi_j(\xi_{lm}^N,\xi_{lm}^N)}=%
        \min\Set[f\in\bignu_j,\ f(\xi_{lm}^N)=1]{\norm{f}}.
        \end{equation*}
    \item We have
        $\Span{\phi_j(\cdot,\xi_{lm}^N)\mid
            \xi_{lm}^N\in X'}=\bignu_j$,
        $j\le j_0$.
    \item The set $\set[\xi_{lm}^N\in X']%
        {\sqrt{\mu_N(\xi_{lm}^N)}\phi_j(\cdot,\xi_{lm}^N)}$
        is a tight frame in $\bignu_j$, $j\le j_0$, i.e.
        for every function $f\in\bignu_j$ we have
        \begin{equation*}\label{label46}
        \norm{f}^2=\summ_{l=1}^N\summ_{m=0}^{2N}\mu_N(\xi_{lm}^N)%
        \left|\innp{f}{\phi_j(\cdot,\xi_{lm}^N)}\right|^2.
        \end{equation*}
    \item
        $$\intt[\omega(\xi)]{_{\setS^2}\phi_j(\xi,\xi_{lm}^N)}=1,$$
        $$\intt[\omega(\xi)]{_{X^2}\phi_j(\xi,\xi_{lm}^N)}=1\qquad
        m_j\le N.$$
    \end{enumerate}}

    {\bf Proof}
    \begin{enumerate}
    \item On the base of \eqref{label35}, every
        $f\in\bignu_j=\Pi_{m_j-1}(\setS^2)$ can be written as
        \begin{equation*}\label{label47}
        f(\eta)=\summ_{n=0}^{m_j-1}\summ_{k=-n}^n
                \innp{f}{Y_{nk}}Y_{nk}(\eta).
        \end{equation*}
        Then
        \begin{multline*}\label{label48}
        \innp{f}{\phi_j(\cdot,\xi_{lm}^N)}=
            \summ_{n=0}^{m_j-1}\summ_{k=-n}^n
            \innp{f}{Y_{nk}}
            \innp{Y_{nk}}{\phi_j(\cdot,\xi_{lm}^N)}= \\
        \summ_{n=0}^{m_j-1}\summ_{k=-n}^n\innp{f}{Y_{nk}}
            \summ_{n'=0}^{m_j-1}\summ_{k'=-n}^n
            \innp{Y_{nk}}{Y_{n'k'}}Y_{n'k'}(\xi_{lm}^N).
        \end{multline*}
        Using the orthonormality property \eqref{label10}
        we obtain that
        \begin{equation*}\label{label49}
            \innp{f}{\phi_j(\cdot,\xi_{lm}^N)}=
                \summ_{n=0}^{m_j-1}\summ_{k=-n}^n
                \innp{f}{Y_{nk}}Y_{nk}(\xi_{lm}^N)=
                f(\xi_{lm}^N).
        \end{equation*}
        Using the discrete orthonormality \eqref{man:3:6} in an
        analogous way we obtain that
        \begin{equation*}\label{label50}
        \innp{f}{\phi_j(\cdot,\xi_{lm}^N)}_X=f(\xi_{lm}^N).
        \end{equation*}
    \item
        \begin{multline*}\label{label51}
        \norm{\phi_j(\cdot,\xi_{lm}^N)}^2=
            \innp{\phi_j(\cdot,\xi_{lm}^N)}
                {\phi_j(\cdot,\xi_{lm}^N)}=
            \phi_j(\xi_{lm}^N,\xi_{lm}^N)=\\
        \summ_{n=0}^{m_j-1}\summ_{k=-n}^n
            \ol{Y_{nk}(\xi_{lm}^N)}Y_{nk}(\xi_{lm}^N)=
            \summ_{n=0}^{m_j-1}(2n+1)\pol_n
                (\xi_{lm}^N\cdot\xi_{lm}^N)=\\
        \summ_{n=0}^{m_j-1}(2n+1)\pol_n(1)=
            \summ_{n=0}^{m_j-1}(2n+1)=m_j^2.
        \end{multline*}
        \item
        Let $f\in\bignu_j$ with $f(\xi_{lm}^N)=1$. Then
        \begin{equation}\label{man:5:16}
        1=\summ_{n=0}^{m_j-1}\summ_{k=-n}^n
        \alpha_{nk}(f)Y_{nk}(\xi_{lm}^N).
        \end{equation}
        Applying the Cauchy-Schwarz inequality we obtain that
        \begin{equation}\label{man:5:17}
        1\le \big(\summ_{n=0}^{m_j-1}\summ_{n=-n}^n|\alpha_{nk}(f)|^2\big)
            \big(\summ_{n=0}^{m_j-1}\summ_{k=-n}^n
            \ol{Y_{nk}(\xi_{lm}^N)}Y_{nk}(\xi_{lm}^N)\big),
        \end{equation}
        with equality attained for $\tilde f$, for which the vectors
        $$\Set{\alpha_{nk}(\tilde f)}_{n=0;k=-n}^{m_j-1;n},\quad
        \Set{Y_{nk}(\xi_{lm}^N)}_{n=0;k=-n}^{m_j-1;n}$$
        satisfy the following conditions: there exists a constant $\alpha\in \setC$ such that $\alpha_{nk}(\tilde f)=
        \alpha\ol{Y_{nk}(\xi_{lm}^N)}$ for $n=\ol{0,m_j-1}$,
        $k=\ol{-n,n}$.

        From { 2.} we deduce that
        $\alpha=\frac{1}{m_j^2}$, and the corresponding $\tilde
        f=\frac{1}{m_j^2}\phi_j(\cdot,\xi_{lm}^N)$.

        For all $f\in\bignu_j$ with
        $f(\xi_{lm}^N)=1$, from \eqref{man:5:17}  and the Parseval's equation it follows
        that $\norm{f}^2\ge\norm{\tilde f}^2=\frac{1}{m_j^2}$.

        Thus
        $$\min\{\norm{f} : \phi\in\bignu_j,\ f(\xi_{lm}^N)=1\}=
        \frac{1}{m_j}=\frac{\norm{\phi_j(\cdot,\xi_{lm}^N})}
            {\phi_j(\xi_{lm}^N,\xi_{lm}^N)}.$$

    \item Writing $f\in\bignu_j$ as its Fourier--Laplace sum
        $$f=\summ_{n=0}^{m_j-1}\summ_{k=-n}^n
            \alpha_{nk}(f)Y_{nk},$$
        on the base of \eqref{man:5:10}, $\alpha_{nk}$ can be
        computed as
        $$\alpha_{nk}(f)=\summ_{k'=1}^{N}\summ_{j=0}^{2N}
        f(\xi_{k'j}^N)\cdot
        \ol{Y_{nk}(\xi_{k'j}^N)}\mu_N(\xi_{k'j}^N).$$
        This means that
        \begin{multline}\label{man:5:18}
        f(\xi)=\summ_{n=0}^{m_j-1}\summ_{k=-n}^n\left(
            \summ_{k'=1}^{N}\summ_{j=0}^{2N}
            f(\xi_{k'j}^N)\ol{Y_{nk}(\xi_{k'j}^N)}
            \mu_N(\xi_{k'j}^N)
            \right)Y_{nk}(\xi)=
        \\
        \summ_{k'=1}^{N}\summ_{j=0}^{2N}f(\xi_{k'j}^N)
            \phi_j(\xi,\xi_{k'j}^N)\cdot\mu_N(\xi_{k'j}^N)=
                \innp{f(\cdot)}{\phi_j(\xi,\cdot)}_X.
        \end{multline}
        Hence $\bignu_j$ is spanned by the functions
        $\phi_j(\cdot,\xi_{k'j}^N)$,
        $\xi_{k'j}^N\in X',\, j\le j_0<N$.
    \item For every function $f\in\bignu_j$ we write
        $$\norm{f}^2=\innp{f}{f}=
            \summ_{n=0}^{m_j-1}\summ_{k=-n}^n
            \ol{\alpha_{nk}(f)}\innp{f}{Y_{nk}}.$$
        From \eqref{man:5:18} it follows that
        $$\norm{f}^2=\summ_{n=0}^{m_j-1}\summ_{k=-n}^n
        \summ_{k'=1}^{N}\summ_{j=0}^{2N}\ol{\alpha_{nk}(f)}
        f(\xi_{k'j}^N)\mu_N(\xi_{k'j}^N)
        \innp{\phi_j(\cdot,\xi_{k'j}^N)}{Y_{nk}}.$$
        Using {1.} we obtain that
        $\innp{\phi_j(\cdot,\xi_{k'j}^N)}{Y_{nk}}=
        \ol{Y_{nk}(\xi_{k'j})}$ and
        $f(\xi_{k'j}^N)=\innp{f}{\phi_j(\cdot,\xi_{k'j}^N)}$,
        then
        \begin{multline*}
        \norm{f}^2=\summ_{n=0}^{m_j-1}\summ_{k=-n}^n
        \summ_{k'=1}^{N}\summ_{j=0}^{2N}
        \ol{\alpha_{nk}(f)}
        \innp{f}{\phi_j(\cdot,\xi_{k'j}^N)}
        \ol{Y_{nk}(\xi_{k'j})}
        \mu_N(\xi_{k'j}^N)=
        \\
        \summ_{k'=1}^{N}\summ_{j=0}^{2N}
        \innp{f}{\phi_j(\cdot,\xi_{k'j}^N)}
        \mu_N(\xi_{k'j}^N)
        \summ_{n=0}^{m_j-1}\summ_{k=-n}^n
        \ol{\alpha_{nk}(f)}
        \ol{Y_{nk}(\xi_{k'j}^N)}=
        \\
        \summ_{k'=1}^{N}\summ_{j=0}^{2N}
        \innp{f}{\phi_j(\cdot,\xi_{k'j}^N)}
        \ol{f(\xi_{k'j}^N)}
        \mu_N(\xi_{k'j}^N)=
        \\
        \summ_{k'=1}^{N}\summ_{j=0}^{2N}
        \innp{f}{\phi_j(\cdot,\xi_{k'j}^N)}
        \ol{\innp{f}{\phi_j(\cdot,\xi_{k'j}^N)}}
        \mu_N(\xi_{k'j}^N)=
        \\
        \summ_{k'=1}^{N}\summ_{j=0}^{2N}
        \left|\innp{f}{\sqrt{\mu_N(\xi_{k'j}^N)}
            \phi_j(\cdot,\xi_{k'j}^N)}\right|^2.
        \end{multline*}
        This means that
        $$\left\{\sqrt{\mu_N(\xi_{lm}^N)}
            \phi_j(\cdot,\xi_{lm}^N),\quad\xi_{lm}^N\in
            X\right\}$$
        is a tight frame in $\bignu_j$ with frame bound $A=1$.
    \item From the orthonormality property \eqref{man:3:6} and
        taking into account that $Y_{00}(\xi)=1$ we have
        \begin{multline*}
        \int_{\setS^2}\phi_j(\xi,\xi_{lm}^N)=
            \intt[\omega(\xi)]{_{\setS^2}
            \summ_{n=0}^{m_j-1}\summ_{k=-n}^n
            \ol{Y_{nk}(\xi_{lm}^N)}Y_{nk}(\xi)}=
            \\
            \summ_{n=0}^{m_j-1}\summ_{k=-n}^n
            \ol{Y_{nk}(\xi_{lm}^N)}
            \intt[\omega(\xi)]{_{\setS^2}}{\ol{Y_{00}(\xi)}Y_{nk}(\xi)}=
            \\
            \summ_{n=0}^{m_j-1}\summ_{k=-n}^n
            \ol{Y_{nk}(\xi_{lm}^N)}\de_{0n}\de_{0k}=1.
        \end{multline*}
        The proof regarding to the discrete integral
        $\int_X\phi_j(\cdot,\xi_{lm}^N)$ is analogous.
    \end{enumerate}
}
\noindent We mention that for the special choice $m_j=j$, with the corresponding frame system
$$\left\{\sqrt{\mu_N(\xi_{lm}^N)}
            \phi_j(\cdot,\xi_{lm}^N),\quad\xi_{lm}^N\in
            X\right\},$$
the analogue of Fej\'er and de la Val\'ee-Poussin summation methods  introduced in \cite{pap} can be expressed.
\section{The harmonic polynomial wavelet space}\label{sec6}
We define the wavelet space $\wlsp_j$, $j\in\setN$, $j\le j_0-1$ as the direct sum
\begin{equation}\label{man:6:1}
\wlsp_j=\dsumm_{n=m_j}^{m_{j+1}-1}\mathcal{H}_n.
\end{equation}
The dimension of $\wlsp_j$ is
$\dim \wlsp_j=\dim\bignu_{j+1}-\dim\bignu_j=m_{j+1}^2-m_j^2$.
Let consider the wavelets $\tilde\psi_j$ for $\wlsp_j$
\begin{equation}\label{man:6:2}
\tilde\psi=\summ_{n=m_j}^{m_{j+1}-1}\summ_{k=-n}^nY_{nk}
\end{equation}
and the ''weighted'' wavelets
\begin{equation}\label{man:6:3}
\psi_j(\cdot,\xi_{lm}^N)=
\summ_{n=m_j}^{m_{j+1}-1}\summ_{k=-n}^n
\ol{Y_{nk}(\xi_{lm}^N)}Y_{nk}\in\wlsp_j,\quad
\xi_{lm}^N\in X'.
\end{equation}
From the addition theorem \eqref{man:2:4} it follows that
\begin{equation}\label{man:6:4}
\psi_j(\cdot,\xi_{lm}^N)=
\summ_{n=m_j}^{m_{j+1}-1}\summ_{k=-n}^n
(2n+1)P_n(\cdot,\xi_{lm}^N),
\end{equation}
consequently $\psi_j(\cdot,\xi_{lm}^N)$, with $\xi_{lm}^N\in X'$ is
real valued. We
summarize the properties of $\psi_j(\cdot,\xi_{lm}^{N})$  in the following theorem.

{\thm{\
      \begin{enumerate}
    \item For every $f\in\wlsp_j$, $j\le j_0-1$ and $\xi_{lm}^N\in X'$ there are  valid the reproducing properties
        $$\innp{f}{\psi_j(\cdot, \xi_{lm}^N)}=f(\xi_{lm}^N),
         $$
        $$\innp{f}{\psi_j(\cdot,
            \xi_{lm}^N)}_X=f(\xi_{lm}^N).\quad
            $$
    \item  The following orthogonality properties hold:
        $$\innp{\phi_j(\cdot, \xi_{lm}^N)}{\psi_j(\cdot,
        \xi_{pq}^N)}=0,\qquad \innp{\phi_j(\cdot, \xi_{lm}^N)}{\psi_j(\cdot,
        \xi_{pq}^N)}_X=0 $$ $$\xi_{lm}^N,\xi_{pq}^N\in X',\quad \xi_{lm}^N \ne \xi_{pq}^N.$$

    \item It holds
        $$\norm{\psi_j(\cdot,\xi_{lm}^N)}=\sqrt{m_{j+1}^2-m_j^2},
        \quad\psi_j(\xi_{lm}^N,\xi_{lm}^N)=m_{j+1}^2-m_j^2.$$
    \item $$\frac{\norm{\psi_j(\cdot,\xi_{lm}^N)}}
            {\psi_j(\xi_{lm}^N,\xi_{lm}^N)}=
            \frac{1}{\sqrt{m_{j+1}^2-m_j^2}}=
            \min{\{\norm{f}:\, f\in\wlsp_j,\,
             f(\xi_{lm}^N)=1\}},$$
    \item $$\Span{\psi_j(\cdot,\xi_{lm}^N),\, \xi_{lm}^N\in X'}=
        \wlsp_j, \ j\le j_0-1, $$
        $$f(\xi)=\innp{f}{\psi_j(\xi,\cdot)}_X,\
        f\in\wlsp_j,\ \xi\in\setS^2,\ j\le j_0<N.$$
    \item The set
        $$\left\{\sqrt{\mu_N(\xi_{lm}^N)}\psi_j(\cdot,
        \xi_{lm}^N):\ \xi_{lm}^N\in X'\right\}$$
        is a tight frame in $\wlsp_j$, $j\le
        j_0$, i.e. for every $f\in\wlsp_j$,
        $$\norm{f}^2=\summ_{l=1}^N\summ_{m=0}^{2N}
        \mu_N(\xi_{lm}^N)|
        \innp{f}{\phi_j(\cdot,\xi_{lm}^N)}|^2.$$
    \item $$\intt[(\xi)]{_{\setS^2}\psi_j(\xi\cdot\xi_{lm}^N)}=0,$$
        $$\intt[\mu_N(\xi)]{_X\psi_j(\xi\cdot\xi_{lm}^N)}=
        0,\ j\le j_0<N, m_j>1.$$
\end{enumerate}}
{\bf Proof}
\begin{enumerate}
    \item Every $f\in\wlsp_j$ can be written as
        $$f=\summ_{n=m_j}^{m_{j+1}-1}\summ_{k=-n}^n
            \alpha_{nk}(f)Y_{nk},\quad
            \alpha_{nk}=\innp{f}{Y_{nk}},$$ hence the scalar product it is equal to:
        \begin{multline*}
        \innp{f}{\psi_j(\cdot,\xi_{lm}^N)}=
        \summ_{n=m_j}^{m_{j+1}-1}\summ_{k=-n}^n
        \innp{f}{Y_{nk}}\innp{Y_{nk}}{\psi_j(\cdot,\xi_{lm}^N)}=
        \\
        \summ_{n=m_j}^{m_{j+1}-1}\summ_{k=-n}^n
        \innp{f}{Y_{nk}}
        \summ_{n'=m_j}^{m_{j+1}-1}\summ_{k'=-n'}^{n'}
        \ub{\innp{Y_{nk}}{Y_{n'k'}}}_{\de_{nn'}\de_{kk'}}
        {Y_{n'k'}(\xi_{lm}^N)}=
        \\
        \summ_{n=m_j}^{m_{j+1}-1}\summ_{k=-n}^n
        \innp{f}{Y_{nk}}
        {Y_{nk}(\xi_{lm}^N)}=f(\xi_{lm}^N).
        \end{multline*}
         For the discrete scalar product the proof is similar.
    \item From the definition of $\phi_j$, $\psi_j$ and the continuous and discrete orthogonality of spherical functions it follows these orthogonality properties. Indeed, for the continuous case we have:
    \begin{multline*}
        \innp{\phi_j(\cdot, \xi_{lm}^N)}
            {\psi_j(\cdot,\xi_{pq}^N)}=
        \\
        \innp{\summ_{n=0}^{m_j-1}\summ_{k=-n}^n
            Y_{nk}(\xi)\ol{Y_{nk}(\xi_{lm}^N)}}
        {\summ_{n'=m_j}^{m_{j+1}-1}\summ_{k'=-n'}^{n'}
            Y_{n'k'}(\xi)\ol{Y_{n'k'}(\xi_{pq}^N)}}=
        \\
        \summ_{n=0}^{m_j-1}\summ_{k=-n}^n
        \summ_{n'=m_j}^{m_{j+1}-1}\summ_{k'=-n'}^{n'}
        \ub{\innp{Y_{nk}(\xi)}{Y_{n'k'}(\xi)}}_{\de_{nn'}\cdot
        \de_{kk'}=0}\ol{Y_{nk}(\xi_{lm}^N)}Y_{n'k'}(\xi_{pq}^N)=0.
        \end{multline*}
        The proof for the discrete scalar product is analogous.
    \item \begin{multline*}
        \norm{\psi_j(\cdot,\xi_{lm}^N)}^2=
        \innp{\psi_j(\cdot,\xi_{lm}^N)}
        {\psi_j(\cdot,\xi_{lm}^N)}=
        \\
        \psi_j(\xi_{lm}^N,\xi_{lm}^N)=
        \summ_{n=m_j}^{m_{j+1}-1}\summ_{k=-n}^n
        Y(\xi_{lm}^N)\ol{Y(\xi_{lm}^N)}=
        \\
        \summ_{n=m_j}^{m_{j+1}-1}(2n+1)
        \ub{P_n(\xi_{lm}^N\cdot\xi_{lm}^N)}_1=
        \summ_{n=m_j}^{m_{j+1}-1}(2n+1)=(m_{j+1})^2-m_j^2.
        \end{multline*}
    \item Let $f\in\wlsp_j$ with $f(\xi_{lm}^N)=1$, then
        \begin{equation}\label{star}
        1=\summ_{n=m_j}^{m_{j+1}-1}\summ_{k=-n}^{n}
        \alpha_{nk}(f)Y_{nk}(\xi_{lm}^N).
        \end{equation}
        Applying the Cauchy-Schwarz inequality we obtain:
        $$1\le\big(\summ_{n=m_j}^{m_{j+1}-1}(\summ_{k=-n}^{n}
        |\alpha_{nk}(f)|^2\big)
        \big(\summ_{n=m_j}^{m_{j+1}-1}\summ_{k=-n}^{n}
        \ol{Y_{nk}(\xi_{lm}^N)}Y_{nk}(\xi_{lm}^N)\big),$$
        with equality for $\tilde f$ with
        $$\alpha_{nk}(\tilde f)=
            \alpha\ol{Y_{nk}(\xi_{lm}^N})\quad
        {{n=m_j,\ldots,m_{j+1}-1}\atop{k=-n,\ldots,n}}$$
        $\alpha\in\setC$ constant. Then
        $$\tilde f=\alpha\summ_{n=m_j}^{m_{j+1}-1}\summ_{k=-n}^{n}
        \ol{Y_{nk}(\xi_{lm}^N)}Y_{nk}(\xi_{lm}^N),$$
        with
        $\alpha=\frac{1}{(m_{j+1}^2-m_j^2)}$.
        If $f\in\wlsp_j$ with $f(\xi_{lm}^N)=1$, then
        $$\norm{f}^2\ge\norm{\tilde f}^2=\frac{1}
        {(m_{j+1}^2-m_j^2)}.$$
        Thus
        $$\min{\left\{\norm{f}:\ f\in\wlsp_j,\,
        f(\xi_{lm}^N)=1\right\}}=
        \frac{1}{\sqrt{(m_{j+1}^2-m_j^2)}}=
        \frac{\norm{\psi_j(\cdot,\xi_{lm}^N)}}
            {\psi_j(\xi_{lm}^N,\xi_{lm}^N)}.$$
    \item Writing $f\in\wlsp_j$, as its Fourier-Laplace sum
        $$f=\summ_{n=m_j}^{m_{j+1}-1}\summ_{k=-n}^{n}
        \alpha_{nk}(f)Y_{nk}
        \text{\quad for \quad }m_{j+1}\le N,$$
        $$\alpha_{nk}(f)=\innp{f}{Y_{nk}}_X=
        \summ_{k'=1}^N\summ_{j=0}^{2N}
        f(\xi_{k'j}^N)\ol{Y_{nk}(\xi_{k'j}^N)}
        \mu_N(\xi_{k'j}^N).$$
        This means that
        \begin{multline}\label{man:6:5}
        f(\xi)=\summ_{n=m_j}^{m_{j+1}-1}\summ_{k=-n}^{n}
        \left(\summ_{k'=1}^N\summ_{j=0}^{2N}
        f(\xi_{k'j}^N)\ol{Y_{nk}(\xi_{k'j}^N)}
        \mu_N(\xi_{k'j}^N)\right)Y_{nk}(\xi)=
        \\
        \summ_{k'=1}^N\summ_{j=0}^{2N}\left(
        \summ_{n=m_j}^{m_{j+1}-1}\summ_{k=-n}^{n}
        \ol{Y_{nk}(\xi_{k'j}^N)}Y_{nk}(\xi)\right)
        f(\xi_{k'j}^N)\mu_N(\xi_{n'j}^N)=
        \int_Xf(\cdot)\psi_j(\xi,\cdot)d\mu_N,
        \end{multline}
        from this it follows that
        $$\Span{\psi_j(\cdot,\xi_{lm}^N):\ \xi_{lm}^N\in X'}=
        \wlsp_j\text{\quad if\quad }
        {m_{j+1}\le N\atop j\le j_0-1}.$$
    \item \begin{multline*}
        \norm{f}^2=\innp{f}{f}=
        \summ_{n=m_j}^{m_{j+1}-1}\summ_{k=-n}^{n}
        \ol{\alpha_{nk}(f)}\innp{f}{Y_{nk}}=
        \\
        \summ_{n=m_j}^{m_{j+1}-1}\summ_{k=-n}^{n}
        \ol{\alpha_{nk}(f)}\innp{
        \summ_{k'=1}^N\summ_{j=0}^{2N}
        f(\xi_{k'j}^N)\phi_j(\cdot,\xi_{k'j}^N)
        \mu_N(\xi_{k'j}^N)}{Y_{nk}}=
        \\
        \summ_{n=m_j}^{m_{j+1}-1}\summ_{k=-n}^{n}
        \ol{\alpha_{nk}(f)}
        \summ_{k'=1}^N\summ_{j=0}^{2N}
        \ub{f(\xi_{k'j}^N)}_{\innp{f}{\phi_j(\cdot,\xi_{k'j}^N)}}
        \ub{{\innp{\phi_j(\cdot,\xi_{k'j}^N)}{Y_{nk}}}}_{
        \ol{Y_{nk}(\xi_{k'j})}}=
        \\
        \summ_{k'=1}^N\summ_{j=0}^{2N}
        \innp{f}{\phi_j(\cdot,\xi_{k'j}^N)}\mu_N(\xi_{k'j}^N)
        \ub{\summ_{n=m_j}^{m_{j+1}-1}\summ_{k=-n}^{n}
        \ol{\alpha_{nk}(f)}
        \ol{Y_{nk}(\xi_{k'j})}}_{\ol{f(\xi_{k'j}^N)}}=
        \\
        \summ_{k'=1}^N\summ_{j=0}^{2N}
        \innp{f}{\phi_j(\cdot,\xi_{k'j}^N)}\mu_N(\xi_{k'j}^N)
        \ol{\innp{f}{\phi_j(\cdot,\xi_{k'j}^N)}}=
        \\
        \summ_{k'=1}^N\summ_{j=0}^{2N}
        \left|\innp{f}{\sqrt{\mu_N(\xi_{k'j}^N)}
            \phi_j(\cdot,\xi_{k'j}^N)}\right|^2.
        \end{multline*}
\end{enumerate}
}

\end{document}